%% file: main.tex
\documentclass[11pt]{article}
\usepackage{graphicx} 

\usepackage[maxbibnames=9]{biblatex}
\usepackage{fancyhdr}
\usepackage{isomath}
\usepackage{mathtools} 
\usepackage{amsbsy}
\usepackage{amssymb}
\usepackage{amscd,amsfonts}
\usepackage{bigints}
\usepackage{graphicx}
\usepackage{verbatim}
\usepackage{euscript}
\usepackage{alltt}
\usepackage{stmaryrd}
\usepackage{relsize}
\usepackage{enumerate}
\usepackage{url}
\usepackage[stable]{footmisc}
\usepackage{breakurl}
\usepackage{hyperref}
\usepackage{comment}
\usepackage[font=small,labelfont=bf]{caption}
\usepackage[font=small,labelfont=bf]{subcaption}
\usepackage{float}
\usepackage{mwe}
\usepackage{algorithm}
\usepackage{algpseudocode}
\usepackage{xcolor}
\usepackage[super]{nth}
\usepackage{soul}

\DeclareGraphicsExtensions{.eps,.pdf}
\input{temp-definitions}

\usepackage[margin=1in]{geometry}

\date{}
\begin{document}
\title{Mid-surface scaling invariance of some bending strain measures}

\author{Amit Acharya\thanks{Department of Civil \& Environmental Engineering, and Center for Nonlinear Analysis, Carnegie Mellon University, Pittsburgh, PA 15213, email: acharyaamit@cmu.edu.}}

\maketitle
\begin{abstract}
\noindent The mid-surface scaling invariance of bending strain measures proposed in \cite{ach_ksb} is discussed in light of the work of \cite{neff}.

\end{abstract}

\section{Introduction}
This brief note discusses the mid-surface scaling invariance of three nonlinear measures of pure bending strain, as introduced in \cite{ach_ksb} and physically motivated therein more than 20 years ago, in light of the recent work of \cite{neff} where the said invariance is introduced. 

It is shown that one of the strain measures introduced in \cite{ach_ksb} possesses scaling invariance, and the other two are easily modified to have the invariance as well. There has been a recent surge of interest in such matters, as can be seen from the works of \cite{neff, virga2024pure, vitral2023dilation}.

We use the notation of \cite{ach_ksb}: a shell mid-surface is thought of as a 2D surface in ambient 3D space (the qualification `mid-surface' will not be used in all instances; it is hoped that the meaning will be clear from the context). Both the
reference and deformed shells are parametrized by the same coordinate system $((\xi^\alpha), \alpha = 1,2)$ (convected
coordinates). Points on the reference geometry are denoted generically by $\bfX$ and on the deformed
geometry by $\bfx$. The reference unit normal is denoted by $\bfN$ and the unit normal on the deformed
geometry by $\bfn$. A subscript comma refers to partial differentiation, e.g. $\frac{\p ()}{\p \xi^\alpha} = ()_{,\alpha}$. Summation over repeated indices will be assumed. The convected coordinate basis vectors in the reference geometry will be referred to by the symbols $(\bfE_\alpha)$ and those in the deformed geometry by $(\bfe_\alpha)$, $\alpha =1,2$, with corresponding dual bases $(\bfE^\alpha), (\bfe^\alpha)$, respectively. A suitable number of dots placed between two tensors represent the operation of contraction, while the symbol $\otimes$ will represent a tensor product. The deformation gradient will be denoted by $\bff = \bfe_\alpha \otimes \bfE^\alpha$ and admits the right polar decomposition $\bff = \bfr \cdot \bfU$, where $\bfU(\bfX): T_\bfX \to T_\bfX$ and $\bfr(\bfX): T_\bfX \to T_\bfx$, where $T_\bfc$ represents the tangent space of the shell at the point $\bfc$. The curvature tensor on the deformed shell is denoted as $\bfb = \bfn_{,\beta} \otimes \bfe^\beta$ and that on the undeformed shell as $\bfB = \bfN_{,\beta} \otimes \bfE^\beta$.

\section{Some measures of pure bending and their invariance under mid-surface scaling}
In \cite{ach_ksb} three measures of bending strain were proposed, given by
\begin{subequations}
    \begin{align}
        \widetilde{\bfK} & = \left(\bfE_\alpha \cdot \bfU  \cdot \bfr^T \cdot \bfn_{, \beta} - \bfE_\alpha \cdot \bfU \cdot \bfN_{,\beta} \right) \bfE^\alpha \otimes \bfE^\beta = \bff^T \cdot \bfb \cdot \bff - \bfU \cdot \bfB \label{eq:unsym_gen_KSB}\\
        \check{\bfK} & = \left( \bfE_\alpha \cdot \bfU  \cdot \bfr^T \cdot \bfn_{, \beta} - \half \left( \bfE_\alpha \cdot \bfU \cdot \bfN_{,\beta} + \bfE_\beta \cdot \bfU \cdot \bfN_{,\alpha}\right) \right) \bfE^\alpha \otimes \bfE^\beta \notag\\
        & = \bff^T \cdot \bfb \cdot \bff - (\bfU \cdot \bfB)_{sym}\label{eq:sym_gen_KSB}\\
        \overline{\bfK} & = \left(\bfE_\alpha \cdot \bfr^T \cdot \bfn_{, \beta} - \bfE_\alpha \cdot \bfN_{,\beta} \right) \bfE^\alpha \otimes \bfE^\beta = \bfr^T \cdot \bfb \cdot \bff - \bfB. \label{eq:r_bend}
    \end{align}
\end{subequations}
Equation \eqref{eq:r_bend} was unnumbered in that work, as the main emphasis was to obtain a nonlinear generalization of the Koiter-Sanders-Budiansky bending strain measure \cite{koiter, b-s}; $\widetilde{\bfK}$ is introduced as Equation (8) and $\check{\bfK}$ as Equation (10) in \cite{ach_ksb}.

In  \cite{neff} a physically natural requirement of invariance of bending strain measure under simple scalings of the form 
$$
\bfx \to a \bfx, \qquad 0 < a \in \R
$$ 
is introduced (for plates, but the requirement is natural for shells as well) and it is shown that the measures $\widetilde{\bfK}, \check{\bfK}$ are not invariant under such a scaling. The measure $\overline{\bfK}$ is not discussed in \cite{neff}.

It is straightforward to see that under the said scaling, the deformation gradient scales as
\[
\frac{\p \bfx}{\p \bfX} = \bfr \cdot \bfU = \bff \qquad \to \qquad a \bff = \bfr \cdot (a \bfU) = a \frac{\p \bfx}{\p \bfX},
\]
resulting in the bending measures scaling as
\[
\widetilde{\bfK} \to a \widetilde{\bfK}; \qquad \check{\bfK} \to a \check{\bfK}; \qquad \overline{\bfK} \to \overline{\bfK}.
\]
Thus, the bending strain measure $\overline{\bfK}$ from \cite{ach_ksb}, not discussed by \cite{neff}, is actually {\bf \textit{invariant under scaling deformations of the deformed shell mid-surface}}. Furthermore, the simple modifications of the measures $\widetilde{\bfK}, \check{\bfK}$ to
\begin{subequations}
    \begin{align}
        \widetilde{\bfK}_{mod} & = \frac{1}{|\bfU|}\left(\bfE_\alpha \cdot \bfU  \cdot \bfr^T \cdot \bfn_{, \beta} - \bfE_\alpha \cdot \bfU \cdot \bfN_{,\beta} \right) \bfE^\alpha \otimes \bfE^\beta \label{eq:unsym_gen_KSB_mod}\\
        & = \left(tr\left(\bff^T \bff\right)\right)^{-\half} \left(\bfx_{,\alpha} \cdot \bfn_{, \beta} - \bfE_\alpha \cdot \bfU \cdot \bfN_{,\beta} \right) \bfE^\alpha \otimes \bfE^\beta \notag\\
        \check{\bfK}_{mod} & = \frac{1}{|\bfU|}\left( \bfE_\alpha \cdot \bfU  \cdot \bfr^T \cdot \bfn_{, \beta} - \half \left( \bfE_\alpha \cdot \bfU \cdot \bfN_{,\beta} + \bfE_\beta \cdot \bfU \cdot \bfN_{,\alpha}\right) \right) \bfE^\alpha \otimes \bfE^\beta \label{eq:sym_gen_KSB_mod}\\
        & = \left(tr\left(\bff^T \bff\right)\right)^{-\half} \left( \bfx_{,\alpha} \cdot \bfn_{, \beta} - \half \left( \bfE_\alpha \cdot \bfU \cdot \bfN_{,\beta} + \bfE_\beta \cdot \bfU \cdot \bfN_{,\alpha}\right) \right) \bfE^\alpha \otimes \bfE^\beta \label{eq:sym_gen_KSB_mod},\notag
    \end{align}
\end{subequations}
where
\[
|\bfU| = \sqrt{\bfU : \bfU} = \sqrt{tr\left(\bff^T \bff\right)},
\]
make them mid-surface scaling invariant.
\printbibliography
\end{document}

%% file: temp-definitions.tex
\usepackage{amsmath}
\usepackage{amsbsy}
\usepackage{amssymb}
\usepackage{amscd}
\usepackage{amsfonts}

\newcommand{\R}{\mathbb R}

\newcommand{\bfb}{{\mathbold b}}
\newcommand{\bfc}{{\mathbold c}}

\newcommand{\bfe}{{\mathbold e}}
\newcommand{\bff}{{\mathbold f}}

\newcommand{\bfn}{{\mathbold n}}

\newcommand{\bfr}{{\mathbold r}}

\newcommand{\bfx}{{\mathbold x}}

\newcommand{\bfB}{{\mathbold B}}

\newcommand{\bfE}{{\mathbold E}}

\newcommand{\bfK}{{\mathbold K}}

\newcommand{\bfN}{{\mathbold N}}

\newcommand{\bfU}{{\mathbold U}}

\newcommand{\bfX}{{\mathbold X}}

\newcommand{\beq}{\begin{equation}}
\newcommand{\eeq}{\end{equation}}
\newcommand{\beqs}{\begin{eqnarray}}
\newcommand{\eeqs}{\end{eqnarray}}
\newcommand{\beql}{\begin{equation} \label}
\newcommand{\half}{\frac{1}{2}}

\newcommand{\p}{\partial}